\documentclass{svjour3}
\smartqed  % flush right qed marks, e.g. at end of proof

\usepackage{graphicx}
\usepackage{float}
\usepackage{mathptmx}
\usepackage{amssymb}

\newcommand{\Eqref}[1]{(\ref{#1})}
\newcommand{\Frac}[2]{\mbox{$\displaystyle\frac{#1}{#2}$}}
\newcommand{\Int}{\displaystyle\int}

\newcommand{\sign}{\mathop{\rm sign}\nolimits}
\newcommand{\Vc}[1]{\mbox{\boldmath$#1$}}

\begin{document}
\title{Dynamics of a particle under the Gravitational Potential of a Massive Annulus: properties and equilibrium description}

%\titlerunning{Dynamics under the potential of a massive annulus}

\author{Eva Tresaco \and Antonio Elipe \and Andr\'es Riaguas}

\institute{Eva Tresaco, Antonio Elipe  \at
              Centro Universitario de la Defensa, 50090 Zaragoza. Spain\\
              Grupo de Mec\'anica Espacial - IUMA. Universidad de Zaragoza, Spain
%             \email{etresaco@unizar.es}
                      %  \\
              \and
           Andr\'es Riaguas \at
           Dept. de Matem\'atica Aplicada, Universidad de Valladolid, 42004 Soria. Spain.}

\date{Received: date / Accepted: date}

\maketitle

\begin{abstract}
This paper studies the main features of the dynamics around a massive annular disk. The first part addresses the difficulties
finding an appropriated expression of the gra\-vi\-tational potential of a massive disk, which will be used
to define the differential equations of motion of our dynamical system. The second part describes the main features of the dynamics with special attention to equilibrium of the system.
\end{abstract}

\section{Introduction}\label{intro}
Outer planets of the solar system and probably many of the extrasolar ones have rings. The rings of some planets of our solar system have been widely studied and have been the object on numerous scientific spatial missions; besides, the Asteroid belt can be roughly approximated by a continuous ring, and its global effect on the orbit of Mars is certainly
not negligible and is worth studying. A similar situation occurs in the study of Coulomb's electrostatic field in the electromagnetism, in galactic dynamics, or gravitational lenses\cite{Kon3,Kon1}. These are a few models that motivated previous works found in literature about the dynamics around a circular ring; the aim of this work is to extend those studies by Maxwell \cite{Maxwell}, Scheeres \cite{Scheeres}, Kalvouridis \cite{Kal} and Elipe and coworkers \cite{ArrEli,ArrEliKalPal,ArrEliPal,BrouckeElipe,EliArrKal}, among others, to the case of an annular disk.

\par
We choose a rather simple model but with rich dynamics to illustrate the relevant structures
in the phase-space and understand their implications. This mathematical model has to be considered as a first approach for further studies of more complex dynamical systems \cite{EliTresRiag}. Our goal is to point out that the conclusions we can draw from
the phase-space structure are generic and therefore of interest in the context of more realistic models. Hence, we consider necessary to start with a knowledge of dynamical aspects in some simple cases, extending the results obtained for the case of a solid ring to an annular ring.

Concerning the annulus problem, some papers have been already published \cite{KroNgSny,LasBli,AlbertiVidal1,Fukus}
in which the potential formula
created by a homogeneous annulus disk is described, those expressions are derived in a similar way than classical books (Kellogg \cite{Kellogg}, MacMilan \cite{MacMilan}) about potential theory used for the computation of the potential created by a circular wire. A different approach is given by Kondratyev \cite{Kon1,Kon2} for computing by analytical methods the potential of  two- and three-dimensional gravitating bodies.

In this paper we consider an homogeneous circular annular disk
of major radius $a$ and minor radius $b$ located on the plane $Oxy$ of a Cartesian coordinate system
and with its geometrical center at the origin of coordinates. For this particular body, the potential function is derived in an appropriated closed form that overcomes some difficulties arisen in numerically computing the elliptic integrals involved in the force function, and will allow us to carry out the computation of the most relevant solutions joined with periodic orbits: the equilibrium positions. Besides, we extend the study to the case of two concentric annuli.

\section{The massive disk and its potential function}\label{se:pot}
Our goal is the analysis of the dynamics of an infinitesimal particle moving under the gravitational field of a massive bidimensional annular disk. The potential or the gravitational force due to  this planar body is obtained by combining the corresponding functions for a massive planar disk. Thus, firstly, we need a convenient  expression for the potential of a massive disk.
This potential can be derived in two ways, as a quadrature involving the potential of solid circular wires of growing radius, or directly from its  definition by computing a double integral. There are similarities and relations in the expressions of the potential and the way of deriving them of three different bodies: ring, disk and annulus, we briefly present all of them in the following subsections.

\subsection{Potential of a massive circular wire and disk}
Let us first remind the potential function created by the gravitational attraction of a solid circular wire under the following assumptions. We consider a homogeneous ring of radius $a$ placed on the $Oxy$-plane of a Cartesian coordinate system
and with its center at the origin of coordinates, assume a total mass $M$ and constant density $\sigma$, that is, $M=2\pi a\sigma$. As it was already stated by Gauss, the potential of a circular ring can be expressed as a complete Elliptic Integral of the first kind \cite{Kellogg}.
\begin{figure}[h]
	\centering\includegraphics[width=3.3in]{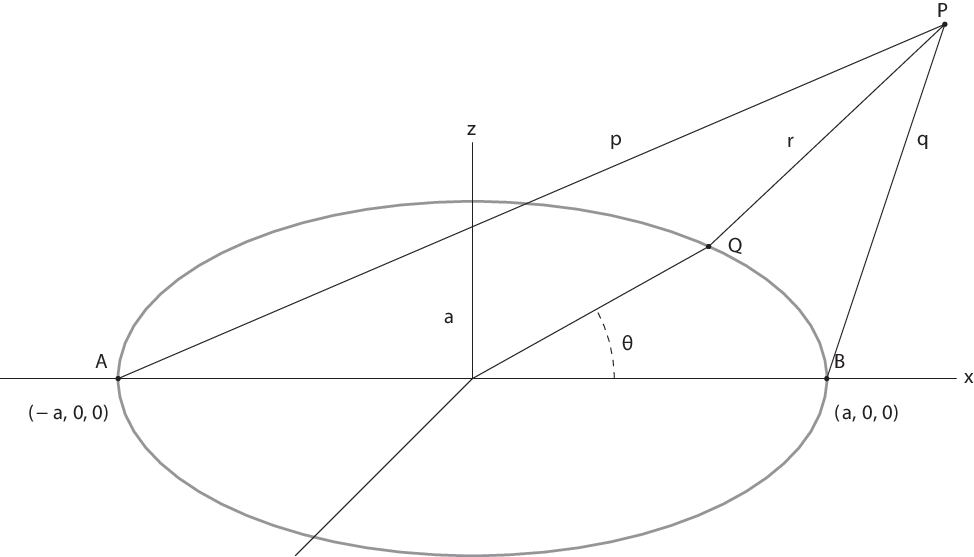}\
	\caption{Wire of radius $a$ acting on a point $P$}
	\label{fig:wire}
\end{figure}

Let $P$ be a mass point of Cartesian coordinates $(x,y,z)$ where the potential is computed  (see Figure~\ref{fig:wire}) and $r$ denotes the distance between $P$ and a differential element of mass $dm$. Let us define the quantities $q$ and $p$, i.e. the shortest and the longest distances between the point P and the ring
\[
  p^2=(x+a)^2+z^2, \qquad q^2=(x-a)^2+z^2,
\]
Introducing in the potential formula angles for the integration, the potential is derived in an elegant form (see Broucke and Elipe  \cite{BrouckeElipe})
\begin{equation}\label{Broucke}
U(P)=-G{\int_0}^{2\pi} \frac{d m}{r}= -\frac{2GM} {\pi p }\int_{0}^{\pi/2}\frac{d\varphi} {\sqrt{1-k^2\sin^2\varphi}} = -\frac{2GM} {\pi p }K(k).
\end{equation}
where $k^2=1-q^2/p^2$, and $K(k)$ is the Complete Elliptic Integral of first kind with modulus $k$. This potential function has a local maximum at the origin, which is unstable (hyperbolic).
\par
Now, let us assume that we have on the $Oxy$ plane a bidimensional disk plate of total mass $M$ and surface density $\sigma$, so $M=\pi\sigma a^2$ and $dm=\sigma ds$.
Let $(r,\theta,z)$ denote the cylindric coordinates for any given point in space $P$ with $r^2=x^2+y^2$ where we want to compute the potential, which is given by
\[
 U(P)=-G\int_D \frac{\sigma d s}d,
\]
where the domain $D$ denotes the disk, and $d$ the distance from a differential mass element $Q$ to the point $P$.

\par
This potential has been already derived by Krough, Ng and Snyder \cite{KroNgSny}  and by Lass and Blitzer \cite{LasBli}. Nevertheless, its closed expression involves elliptic integrals, and the expressions given for the potential cannot be evaluated at significant regions of the space where the potential is a well defined function, or in such a way that may produce wrong evaluations when is numerically computed.
\par
Let us summarize here the approach derived to develop the formulation of the potential function, and how to overcome the difficulties in computing it.
\par
\begin{figure}[h]
	\centering\includegraphics[width=3.5in]{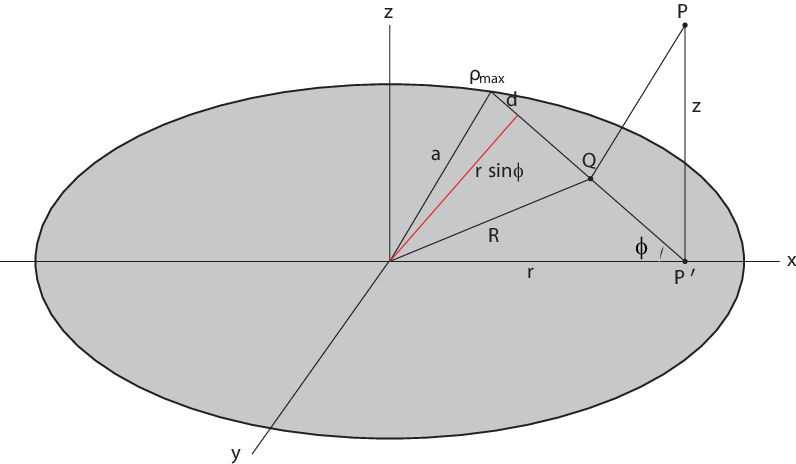}\
	\caption{Disk of radius $a$ acting on a point $P$ with $r\leq a$}
	\label{fig:pot_disk_in}
\end{figure}

Let $P$ be a mass point where the potential is computed, and $Q$ an element of differential mass of the disk at distance $d$ from $P$ (see Fig.~\ref{fig:pot_disk_in}).

We determine first the potential at points $P$ with  $r\leq a$. The integrand in the expression of the potential is in terms of the angle $\phi$ and the distance
$\rho$ from the projection $P'$ of $P$ into the plane $Oxy$ to the point $Q$
\[
    U(P)= -\frac{GM}{\pi a^2}\int_{0}^{2\pi}\int_{0}^{\rho_{max}}\frac{\rho d\rho d\phi}{\|\overline{PQ}\|},
\]
where
$\rho_{max}$ is the maximum distance from $P'$ to $Q$, i.e. when $Q$ is at the boundary of the disk. The distance $\rho_{max}$ and $d$ are given by the formulas
\[\rho_{max}=r\cos\phi+d, \quad d^2=a^2-r^2\sin^2\phi.\]
Hence, the integrand is now
\[
\begin{array}{l}\displaystyle
    U(P)= -\frac{2GM}{\pi a^2}\int_{0}^{\pi}\int_{0}^{\rho_{max}}\frac{\rho d\rho d\phi}{\sqrt{\rho^2+z^2}}\\[2ex]
    \qquad\qquad \displaystyle =\frac{2GM}{\pi a^2}\left(\pi|z|-\int_{0}^{\pi}\sqrt{\left(r\cos\phi+\sqrt{a^2-r^2\sin^2\phi}\right)^2+z^2} d\phi\right).
    \qquad \qquad \end{array}
\]
To solve the quadrature we need to manipulate the expression, we show here the main steps.
\par
 First of all we change the integration variable $\varphi=r\cos\phi+\sqrt{a^2-r^2\sin^2\phi}$; after substitution and simplification the potential yields
\begin{equation}\label{eq:fla_compara}
    U(P)= \frac{2GM}{\pi a^2}\left(\pi|z|-\int_{a-r}^{a+r}\frac{\sqrt{z^2+\varphi^2}(\varphi^2-(r^2-a^2))}{\varphi\sqrt{4r^2a^2-(\varphi^2-(r^2+a^2))^2}}d\varphi\right).
\end{equation}
Next we define a new integration variable $\cos2\psi=(\varphi^2-(r^2+a^2))/2r a$, in order to eliminate the root in the denominator. Thus the potential expression can be written as
\begin{equation}\label{eq:fla_add}
\begin{array}{l}
\displaystyle
U(P)= \frac{2GM}{\pi a^2}\left(\pi|z|-\sqrt{z^2+(a+r)^2}E(k)
\right.\\[3ex]
\qquad\qquad
\displaystyle \left.
+ \sqrt{z^2+(a+r)^2}\frac{(a^2-r^2)}{(a+r)^2}\int_{0}^{\pi/2}\frac{\sqrt{1-k^2\sin^2\psi}}{
1-n^2\sin^2\psi}d\psi
\right),
\end{array}
\end{equation}
where we have introduced the auxiliary quantities
\[k^2=\frac{4ar}{z^2+(a+r)^2},\quad n^2=\frac{4ar}{(a+r)^2}.
\]
Finally, to solve the last integral we multiply both numerator and denominator by the term $\sqrt{1-k^2\sin^2\psi},$ in order to get rid of the root in the denominator.
\par
Lastly, we apply the definitions of the Elliptic Integrals to eventually get the following formula
\begin{equation}\label{eq:pot_disk}
\begin{array}{l}
\displaystyle
    U(P)= \frac{2GM}{\pi a^2}\left(\pi|z|-\sqrt{z^2+(a+r)^2}E(k)- \frac{(a^2-r^2)}{\sqrt{z^2+(a+r)^2}}K(k)
\right.\\[3ex]
\qquad\qquad
\displaystyle \left.
    -\frac{(a-r)}{(a+r)}\frac{z^2}{\sqrt{z^2+(a+r)^2}}\Pi(k,n)\right),
    \end{array}
\end{equation}
with $E(k)$, $K(k)$ and $\Pi(n,k)$ the Complete Elliptic
Integrals of first, second and third kind respectively. Let us remind that this expression is only valid when we compute the potential at points $P$ such that $r\leq a$.
\par
\begin{figure}[h]
	\centering\includegraphics[width=3.5in]{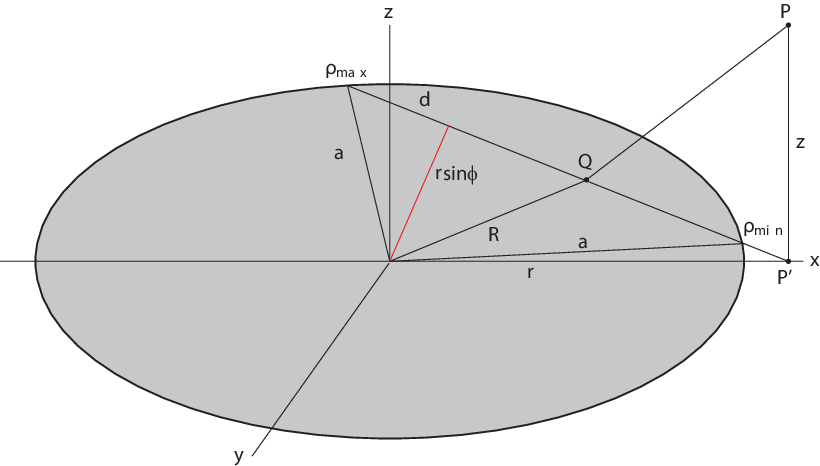}\
	\caption{Disk of radius $a$ acting on a point $P$ with $r>a$}
	\label{fig:pot_disk_out}
\end{figure}

For the case $r>a$, the procedure is similar, see Fig.~\ref{fig:pot_disk_out}.

\noindent
Taking into account that now  $\phi\in[-\arcsin(a/r),\arcsin(a/r)]$ and $\rho\in[\rho_{min},\rho_{max}]$, where
\[
\rho_{max}=r\cos\phi+d, \quad \rho_{min}=r\cos\phi-d, \quad d^2=a^2-r^2\sin^2\phi,
\]
Eq. (\ref{eq:pot_disk}) becomes
\begin{equation}\label{pot_disk9}
\begin{array}{rcl}
U(P)&=& -\Frac{GM}{\pi a^2}\Int_{-\arcsin(a/r)}^{\arcsin(a/r)}\Int_{\rho_{min}}^{\rho_{max}}\Frac{\rho d\rho d\phi}{\sqrt{\rho^2+z^2}}=-\Frac{2GM}{\pi a^2}\Int_{0}^{\arcsin(a/r)}\Int_{\rho_{min}}^{\rho_{max}}\Frac{\rho d\rho}{\sqrt{\rho^2+z^2}}
\\[2ex]
&=& - \Frac{2GM}{\pi a^2}\Big(\Int_{0}^{\arcsin(a/r)}\sqrt{{\rho_{max}}^2+z^2}d\phi-\Int_{0}^{\arcsin(a/r)}\sqrt{{\rho_{min}}^2+z^2}d\phi\Big).
\end{array}
\end{equation}

By means of two changes of variables, one for each integral, $\varphi_1=r\cos\phi+\sqrt{a^2-r^2\sin^2\phi}$, and $\varphi_2=r\cos\phi-\sqrt{a^2-r^2\sin^2\phi}$, the potential function can we written as
\[
\begin{array}{l}
\displaystyle U(P)= -\frac{2GM}{\pi a^2}\Bigg(\int_{r+a}^{\sqrt{r^2-a^2}}\frac{\sqrt{z^2+{\varphi_1}^2}(-{\varphi_1}^2+r^2-a^2)d\varphi_1}{\varphi_1\sqrt{4r^2{\varphi_1}^2-({\varphi_1}^2+r^2-a^2)^2}} \\[3ex]
\quad \displaystyle -\int_{r-a}^{\sqrt{r^2-a^2}}\frac{\sqrt{z^2+{\varphi_2}^2}(-{\varphi_2}^2+r^2-a^2)d\varphi_2}{\varphi_2\sqrt{4r^2{\varphi_2}^2-({\varphi_2}^2+r^2-a^2)^2}}\Bigg)
=-\frac{2GM}{\pi a^2}\int_{r+a}^{r-a}\frac{\sqrt{z^2+\varphi^2}(-\varphi^2+r^2-a^2)d\varphi}{\varphi\sqrt{4r^2\varphi^2-(\varphi^2+r^2-a^2)^2}}.
\end{array}
\]
Now, if we apply the change of variable $\cos2\psi=(\varphi^2-(r^2+a^2))/2r a$ we obtain the expression given in Eq.~\Eqref{eq:fla_add},  so we can conclude that the potential is given by
\begin{equation}\label{eq:pot_disk2}
\begin{array}{l}
\displaystyle
U(P)=\frac{2GM}{\pi a^2}\left(-\sqrt{z^2+(a+r)^2}E(k)- \frac{(a^2-r^2)}{\sqrt{z^2+(a+r)^2}}K(k)\right.
\\[2ex]
\qquad\qquad\displaystyle
\left.
-\frac{(a-r)}{(a+r)}\frac{z^2}{\sqrt{z^2+(a+r)^2}}\Pi(k,n)\right).
\end{array}
\end{equation}
Therefore we obtain the same expression for the potential that we had for the previous case $r\leq a$, see Eq.~\Eqref{eq:pot_disk}, without the term in $|z|$. Both formulas, Eq.~\Eqref{eq:pot_disk} and Eq.~\Eqref{eq:pot_disk2}, can be written as a single one
\begin{equation}
\begin{array}{l}
\displaystyle
    U(P) = \frac{2GM}{\pi a^2}\left(|z|\frac{\pi}{2}(1+\sign(a-r))-\sqrt{z^2+(a+r)^2}E(k)\right.
\\[2.5ex]
\qquad\displaystyle
\left.
- \frac{(a^2-r^2)}{\sqrt{z^2+(a+r)^2}}K(k^2)-\frac{a-r}{(a+r)}\frac{z^2}{\sqrt{z^2+(a+r)^2}}\Pi(k,n) \right).
\end{array}\\[2.5ex]
\end{equation}
Consequently the potential created by an homogeneous bidimensional disk is given by the closed form expression
\begin{equation}\label{eq:ERT-pot1}
U(P)= \frac{2GM}{\pi a^2}\left(|z|\frac{\pi}{2}(1+\sign(a-r))-pE(k)  -\frac{a^2-r^2}{p}K(k)-\frac{(a-r)}{(a+r)}\frac{z^2}{p}\Pi(k,n) \right).
\end{equation}
Nevertheless, this formula does not represent the potential function at every point in the space for which the potential function has a real finite value. Albeit they are valid for the cases $r<a$ and $r>a$, they fail at $r=a$ and $z\ne 0$, because in this case $n=1$ and the  Elliptic Integral of the third kind is not bounded for these values, rendering those formulas useless for the analysis of the dynamics. It is necessary to rewrite this analytic expressions to reduce the loss of information in its computation (Fukushima \cite{Fukus}).

\subsubsection{Proper evaluation through potential reformulation}
In order to circumvent the problems arisen in the computation of the potential expression Eq.~\Eqref{eq:ERT-pot1}, we use several formulas from the Byrd and
Friedman textbook \cite{ByrFri}, and the computational approach by Burlirsch \cite{bul}, Carlson \cite{carl} and Fukushima \cite{Fukus1,Fukus2} to overcome difficulties in evaluating the Elliptic Integrals.

Let us split the potential expression given in Eq.~\Eqref{eq:ERT-pot1} in two parts
\begin{equation}\label{eq:ERT-pot}
 U=\displaystyle\frac{2GM}{\pi a^2}\left(U_1+U_2\right),
\end{equation}
where $U_1 = |z|\displaystyle\frac{\pi}{2}(1+\sign(a-r))-pE(k)-\displaystyle\frac{a^2-r^2}{p}K(k), U_2 = -\displaystyle\frac{(a-r)z^2}{(a+r)p}\Pi(n,k)$\\
We now transform the term $U_2$, which contains the Elliptic Integral of the third kind, by means of the following relations
\[ 1-n^2 =\displaystyle
\frac{(a-r)^2}{(a+r)^2}, \qquad
\displaystyle\frac kn = \displaystyle\frac {a+r}p,
\]
and the relation between the Elliptic Integral of third kind and the Heuman's lambda function $\Lambda_0$ (see [413.01] in Byrd and Friedman \cite{ByrFri})
\[
\Pi(n,k)= \displaystyle \frac{n\pi
\Lambda_0(\phi,n)}{2\sqrt{(n^2-k^2)(1-n^2)}},\qquad \phi=\arcsin\sqrt{\displaystyle\frac{n^2-k^2}{n^2(1-k^2)}}=\arcsin\displaystyle\frac{|
z   |}{q},
\]
to rewrite $U_2$ as
\[
U_2=-\frac{z^2(a-r)}{p(a+r)}\displaystyle \frac{n\pi}
{2\sqrt{(n^2-k^2)(1-n^2)}}\Lambda_0(\phi,n) = -|z|\frac{\pi}{2}\sign(a-r)\Lambda_0(\phi, k).
\]
Finally, it is possible to replace the Heuman's lambda function by a
combination of Elliptic Integrals
\[\Lambda_0(\phi, k)=\frac{2}{\pi}\left(\frac{}{}E(k)F(\phi,k^\prime)+K(k)E(\phi,k^\prime)-K(k)F(\phi,k^{\prime})\right),\] where $k^\prime=\sqrt{1-k^2}$, and
$F(\phi,k^\prime)$ and $E(\phi,k^{\prime})$ are the Incomplete Elliptic Integrals of the first and the second kind,
respectively.\\ Thus $U_2$ can be reformulated as
\[U_2=-|z|\sign(a-r)\left(E(k)F(\phi,k^\prime)+K(k)E(\phi,k^\prime)-K(k)F(\phi,k^{\prime})\right),\]
and the potential function yields
\begin{eqnarray}\label{eq:ERT-pot_def}
U&=&\frac{2GM}{\pi a^2}\left(-pE(k)-\frac{a^2-r^2}{p}K(k)+|z|\left(\frac{\pi}{2}+\frac{\pi}{2}\sign(a-r)\right)\right. \\[1.6ex]
\nonumber &&\left.-|z|\sign(a-r)\Big(E(k)F(\phi,k^\prime)+K(k)E(\phi,k^\prime)-\frac{}{}K(k)F(\phi,k^{\prime})\Big)\right).
\end{eqnarray}
Under this form, the potential function and the force function derived from it can be properly evaluated
at any point in the space where it is defined.

It is worth to remark that a similar gravitating spatial potential of the homogeneous circular disk was obtained by Kondratyev in his two monographs \cite{Kon1,Kon2}, where the author deduces  integral formulas for the gravitational potential of two-dimensional, and cylindrical bodies.

\subsection{Potential of a massive annular disk}

Once we have derived a proper expression for the potential function of the circular plate it is immediate to compute the potential due to an annular disk by subtraction. The potential created by an annulus of radii $a$ and $b$ ($b<a$) is computed from the disk potential by subtracting two concentric disks of radius $a$ and $b$ respectively,
\begin{equation}\label{eq:ERT-pot_ann}
 U(x,y,z;a,b)=U(x,y,z;a)-U(x,y,z;b),
\end{equation}
where $U(x,y,z;a)$ and $U(x,y,z;b)$ are the potentials created by two circular concentric plates of radii $a$ and $b$ (see Eq.~\Eqref{eq:ERT-pot_def}). Note that now
\[
G\sigma=GM/(\pi(a^2-b^2)) = \mu/(\pi(a^2-b^2)).
\]
Hence the potential created by a planar annulus is given by the following formula
\begin{equation}\label{pot_ann}
\begin{array}{rcl}
U=&&\Frac{2\mu}{\pi(a^2-b^2)}\Big(
 -p_aE(k_a)-\Frac{a^2-r^2}{p_a}K(k_a)+|z|\left(\Frac{\pi}{2}+\Frac{\pi}{2}\sign(a-r)\right)\\[2ex]
&&-|z|\sign(a-r)\big[E(k_a)F(\phi_a,k_a^\prime)+K(k_a)E(\phi_a,k_a^\prime)-K(k_a)F(\phi_a,k_a^{\prime}) \big]\\[2ex]
&&+p_bE(k_b)+\Frac{b^2-r^2}{p_b}K(k_b)-|z|\left(\Frac{\pi}{2}+\Frac{\pi}{2}\sign(b-r)\right)\\[2ex]
&&+ |z|\sign(b-r)\big[E(k_b)F(\phi_b,k_b^\prime)+K(k_b)E(\phi_b,k_b^\prime)-K(k_b)F(\phi_b,k_b^{\prime})\Frac{}{}\big]
\Big),
\end{array}
\end{equation}
where we have introduced the following auxiliaries quantities
\[
\begin{array}{lll}
r^2=x^2+y^2, & R^2=x^2+y^2+z^2, \\[2ex]
p_a^2=(a+r)^2+z^2,\quad  & q_a^2=(a-r)^2+z^2, \quad & \\[2ex]
k_a^2=4ar/p_a^2,  &k_a^\prime=\sqrt{1-k_a^2}, &  \phi_a=\arcsin {|
z|}/{q_a},\\[2ex]
p_b^2=(b+r)^2+z^2,  & q_b^2=(b-r)^2+z^2, & \\[2ex]
k_b^2=4br/p_b^2,  &k_b^\prime=\sqrt{1-k_b^2}, & \phi_b=\arcsin {|
z|}/{q_b}.
\end{array}
\]

\section{Dynamics around a circular annulus}\label{se:dyn}
Once we have a convenient expression of the potential we proceed to study the dynamics of an infinitesimal particle moving under the gravitational field of a massive bidimensional circular annulus.
\par
The Lagrangian describing the motion of a particle in space under the attraction of an annular disk is given by
\[
{\cal L}= T - U = \frac12(\dot x^2+\dot y^2 + \dot z^2) - U(x,y,z),
\]
where $U(x,y,z)$ is the potential (Eq.~\Eqref{pot_ann}), and T the kinetic energy. The Euler-Lagrange equations that correspond to this Lagrangian are given by
\begin{equation}\label{eqMot}
\ddot{x} = -U_x, \qquad
\ddot{y} = -U_y, \qquad
\ddot{z} = -U_z,
\end{equation}
where $U_x,U_y,U_z$ denote the partial derivatives $\Big(\partial U/\partial x, \partial U/\partial y, \partial U/\partial z\Big):$
\begin{eqnarray}\label{eq_mvt_ann}
  \nonumber \frac{\partial U}{\partial x}&=&\frac{2\mu}{\pi(a^2-b^2)}\frac{x}{r^2}\left(  \sqrt{R^2+a^2+2ar}[(1-\frac{1}{2}k_a)K(k_a)-E(k_a)]-\right.\\[1ex]
  \nonumber&&\left.\sqrt{R^2+b^2+2br}[(1-\frac{1}{2}k_b)K(k_b)-E(k_b)]\right),\\[1.6ex]
  \frac{\partial U}{\partial y}
\nonumber&=&\frac{2\mu}{\pi(a^2-b^2)}\frac{y}{r^2}\left(  \sqrt{R^2+a^2+2ar}[(1-\frac{1}{2}k_a)K(k_a)-E(k_a)]-\right.\\[1ex]
  \nonumber&&\left.\sqrt{R^2+b^2+2br}[(1-\frac{1}{2}k_b)K(k_b)-E(k_b)]\right),\\[1.6ex]
  \nonumber \frac{\partial U}{\partial z}
&=&-\frac{\mu}{\pi(a^2-b^2)}\left(\frac{2z}{\sqrt{R^2+a^2+2ar}}K(k_a)-2\sign(z)\left(\frac{\pi}{2}+\frac{\pi}{2}\sign(a-r)\right.\right.\\[1.6ex]
\nonumber&&\left.-\sign(a-r)[\frac{}{}(E(k_a)-K(k_a))F(\phi,k_a^{\prime})+K(k_a)E(\phi,k_a^\prime)]\right)-\frac{2z}{\sqrt{R^2+b^2+2br}}K(k_b)\\[1ex]
\nonumber&&+2\sign(z)\left(\frac{\pi}{2}+\frac{\pi}{2}\sign(b-r)-\sign(b-r)\Big[\frac{}{}(E(k_b)-K(k_b))F(\phi_b,k_b^{\prime})\right.\\[1ex]
&&\left.+K(k_b)E(\phi_b,k_b^\prime)\Big]\right)\Bigg).
\end{eqnarray}

The problem is autonomous and hence, the Hamiltonian (hence, the energy) is constant along the motion,
\[
    E= \frac12({\dot x}^2 +{\dot y}^2 +{\dot z}^2) + U(x,y,z) .
\]

\begin{figure}[h]
	\centering\includegraphics[width=2.2in]{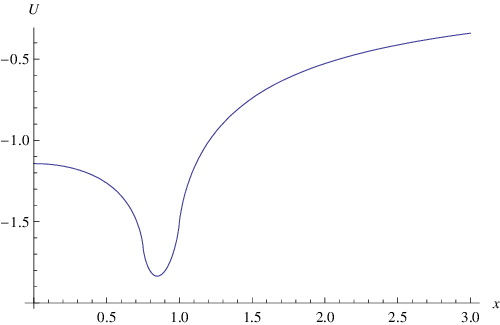}\quad
    \centering\includegraphics[width=2.2in]{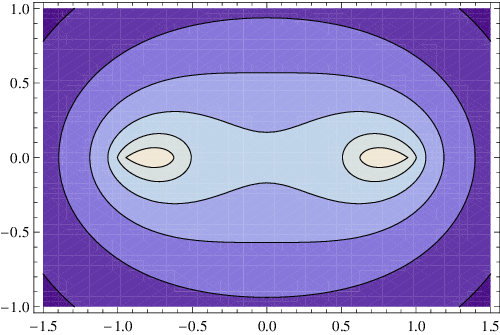}
	\caption{Left: Potential function $U = U(x,0,0)$. Right: Equipotential curves for the function $U(x,y,0)$}
	\label{fig:pot_contours}
\end{figure}

Besides, the potential of the annular disk is symmetric with respect to the three axes $Ox, Oy$ and $Oz$ because of the cylindric symmetry of the problem.
Figure~\ref{fig:pot_contours} (left) illustrates the potential function $U$ along the $Ox$-axis (for the section $y=z=0$), and the right plot depicts equipotential curves of the function $U(x,y,0)$.

It can be observed that the potential function goes to $0$ as $r$ goes to infinity and also the existence of a local maximum at the origin. Let us prove that the origin is a linear unstable equilibrium indeed.
\par
Let us remind that equilibrium points are the simplest invariant objects along with the periodic orbits; they are important not only for their existence, also because they structure the global dynamics of the system.

The equilibrium positions of the system, $\Vc{x}_0$, are the critical points of the potential function given by $\Delta U(\Vc{x}_0)=0.$  The stability of a solution $\Vc{x}_0$ of a Hamiltonian system is given by the variational equations that determine the motion in a neighborhood of the equilibrium position. The equilibrium is spectrally stable if all its eigenvalues are pure imaginary.

In order to prove that the origin is an equilibrium of the system, second derivatives must be computed. We have used their integral expressions instead of a straight derivation of the equation of motion, see Eq.~\Eqref{eq_mvt_ann}, which gives formulas difficult to handle analytically. For that reason we consider now an arbitrary point of the annulus whose polar coordinates can be written as $(\rho\cos\theta,\rho\sin\theta,0)$. The potential created by the annulus on a point $P(x,y,z)$ in space is
\[
U=-G\sigma\int_{b}^{a}\int_{0}^{2\pi}\frac{\rho d\theta d\rho}{{(x^2+y^2+z^2+\rho^2-2\rho x\cos\theta-2\rho y\sin\theta)}^{1/2}}.
\]
Thus, its derivative with respect to $x$ is
\[
U_x= G\sigma\int_{b}^{a}\int_{0}^{2\pi}\frac{\rho(x-\rho\cos\theta)d\theta d\rho}{{(x^2+y^2+z^2+\rho^2-2\rho x\cos\theta-2\rho y\sin\theta)}^{3/2}};
\]
and analogously for $U_y$ and $U_z$. It is clear from those expression that the origin is an equilibrium point of the problem. These expressions are easy to derive again and to particularize them at the origin to get the monodromy matrix whose eigenvalues are (disregarding multiplicity)
\[
 \pm i\sqrt{\pi G \sigma\frac{b-a}{ab}},\quad -\sqrt{2}\sqrt{\pi G\sigma\frac{b-a}{ab}},\quad \sqrt{2}\sqrt{\pi G\sigma\frac{b-a}{ab}}.
 \]
Since there are real eigenvalues, we can conclude that the origin is spectrally unstable.

On the other hand, since the annulus model has axial symmetry, it is natural to use cylindrical coordinates $(r,\lambda,z)$ to have the Lagrangian
\[
{\cal L} = \frac{1}{2} (\dot{r}^2 + r^2\dot{\lambda}^2+\dot{z}^2 ) - U(r,z);
\]
since in this case the angle $\lambda$ is a cyclic variable, its conjugate moment $\Lambda=\partial{\cal L}/\partial \dot{\lambda } = r^2\dot{\lambda}$ is constant and the equations of motion in cylindric coordinates are
\begin{equation}\label{eqMotion}
\begin{array}{l}
\ddot{r}  =   -\partial U /\partial r + \Lambda^2/r^3,\\
\ddot{z}  = -\partial U /\partial z . \\
\end{array}
\end{equation}

In order to find the stationary solutions we have to analyze the equations of motion; due to the complexity of the expressions containing Elliptic Integrals, analytic solutions are not feasible in general. This is the reason why we have studied the equilibrium when the motion is reduced to the $xy$-plane and when the motion is confined to the $Oz$-axis.

 Note that for the following computations we chose the following values for the parameters: $a=1$, $b=0.75$ and $\mu=1$. We proved that we can assume without loss of generality that both the outer radius of the annulus ($a$), and the gravitational constant ($\mu$) are equal to one, and concerning the inner radius of the annulus, we saw that slightly variations of this parameter do not modify qualitatively the nature of the results presented here.

\subsection{Dynamics on the $Oz$-axis}\label{se:dyn_z}

\begin{figure}[htb]
	\includegraphics[width=2.5in]{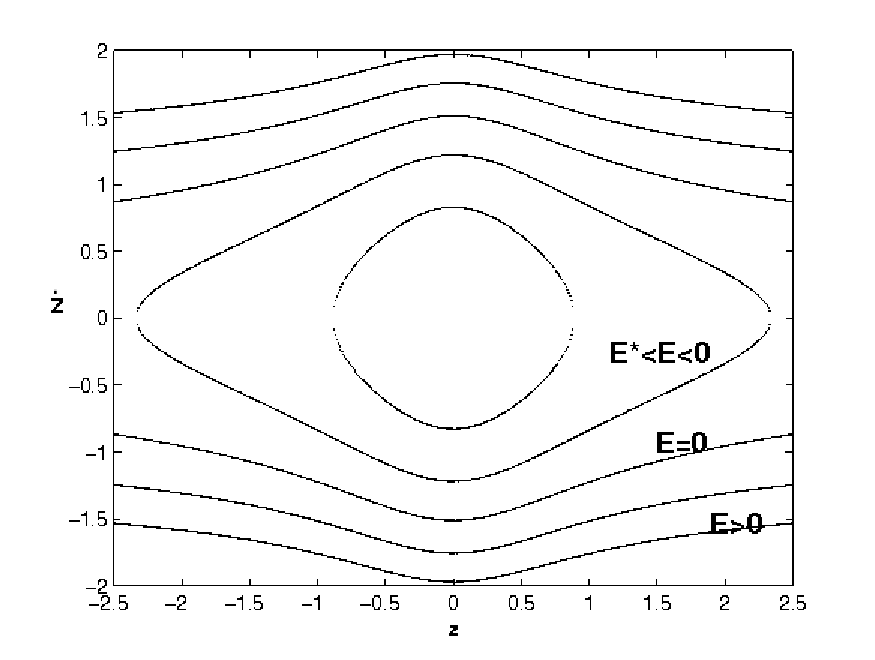}
	\caption{Phase portrait over the $Oz$ axis.}
	\label{fig:phase_portr_zZ}
\end{figure}

Following Kellogg \cite{Kellogg}, the potential of the annulus is a \textsl{single layer potential} with essential discontinuities at
the boundary of the circular plate, but otherwise it is a continuous
function.
Its gradient is also a continuous function everywhere except
at points on the annulus plate. It is not defined for points at
the boundary  and it has a step discontinuity at
points on the plate but outside its boundary. The derivative along the normal direction to the annulus of  $U$  shows a step discontinuity on the annulus
\[
\left(\frac{\partial U}{\partial n}\right)_+ -
\left(\frac{\partial U}{\partial n}\right)_-=4\pi G \sigma.
\]
This discontinuity governs the movement in the $Oz$-axis. Now $r=0$ and the movement is determined by the following differential equation
\[
 \ddot{z}=\frac{\mu z}{\pi(a^2-b^2)}\left(\frac{1}{\sqrt{z^2+a^2}}-\frac{1}{\sqrt{z^2+b^2}}\right).
 \]
Note that it is equal to zero if $z=0$ or $a=b$, but this last case corresponds to the wire problem which has been already study by Broucke and Elipe \cite{BrouckeElipe}. Therefore, the only equilibrium point is $z=0$, because the motion reduces to the $0z$-axis, this corresponds to the origin. The energy at this point is $E^*=-2\mu/(a+b)$, and only for energy values  $E^*<E<0$ we will find periodic orbits. This equilibrium point can be also observed in the phase portrait (Figure~\ref{fig:phase_portr_zZ}).
\subsection{Dynamics on the $xy$--plane}\label{se:dyn_xy}
In this subsection we analyze the equilibrium in the planar case $z=0$, where the effective potential is now a radial function, therefore we are in presence of an integrable problem thanks to the energy and the angular momentum integrals.

Since $z=0$, the potential \Eqref{pot_ann} reduces to
\begin{equation}\label{eq:potefe}
U(r)=\frac{2\mu}{\pi(a^2-b^2)}\Big((b+r)E(k_b)+(b-r)K(k_b)
-(a+r)E(k_a)-(a-r)K(k_a)\Big).
\end{equation}

By using the expressions of the partial derivatives \Eqref{eq_mvt_ann}, we easily compute
\begin{equation}\label{eq:cilEqMot_r_z0}
\begin{array}{l}
\Frac{\partial U}{\partial r}
  = \Frac{2\mu}{\pi(a^2-b^2)} \Frac{1}{r}
  \left(\left[\Frac{(a^2+r^2)}{(a+r)} K(k_a) - (a+r) E(k_a)\right]\right. \\[2.5ex]
\qquad\qquad
- \left.\left[\Frac{(b^2+r^2)}{(b+r)}K(k_b)  -(b+r)E(k_b) \right] \right) .
\end{array}
\end{equation}

Since the system is conservative,
\[
E = T + U(r)=\frac{1}{2}\left({\dot r}^2+\frac{\Lambda^2}{ r^2}
\right) + U(r) = \frac{1}{2}({\dot r}^2+W(r)),
\]
where $W(r) = \Lambda^2/r^2 + U(r)$ is the so-called effective potential. Hence,
\[
\dot r = \sqrt{2 (E - W(r) )}.
\]
Stationary points will be circular solutions of the complete problem $(r,\lambda,z)$, i.e., the critical points of the effective potential, named $r_0$, for values of the angular momentum $\Lambda\neq 0$, will correspond to circular orbits on the plane $Oxy$: $\left(r_0\cos(\Lambda t/{r_0}^2),r_0\sin(\Lambda t/{r_0}^2),0\right)$ of period $T=2\pi{r_0}^2/\Lambda$.

 We are going to find the values of the angular momentum for which periodic orbits e\-xist\-. With this goal in mind, we perform a similar analysis to the one developed by Alberti and Vidal \cite{AlbertiVidal1} to find the equilibrium positions on the equatorial plane of the annulus. One of the motivation of our study is to extend the dynamical approach performed by in the mentioned work to a more complex situation as result of a composition of annulus.

\begin{figure}[h]
\centering
\includegraphics[width= 2.6in]{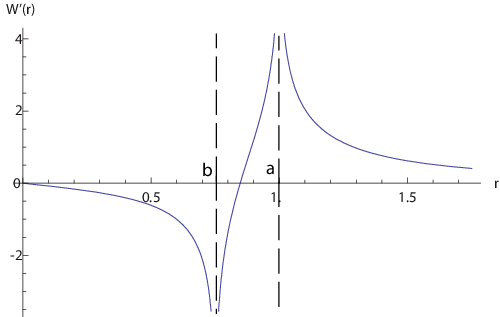}
\includegraphics[width= 2.1in]{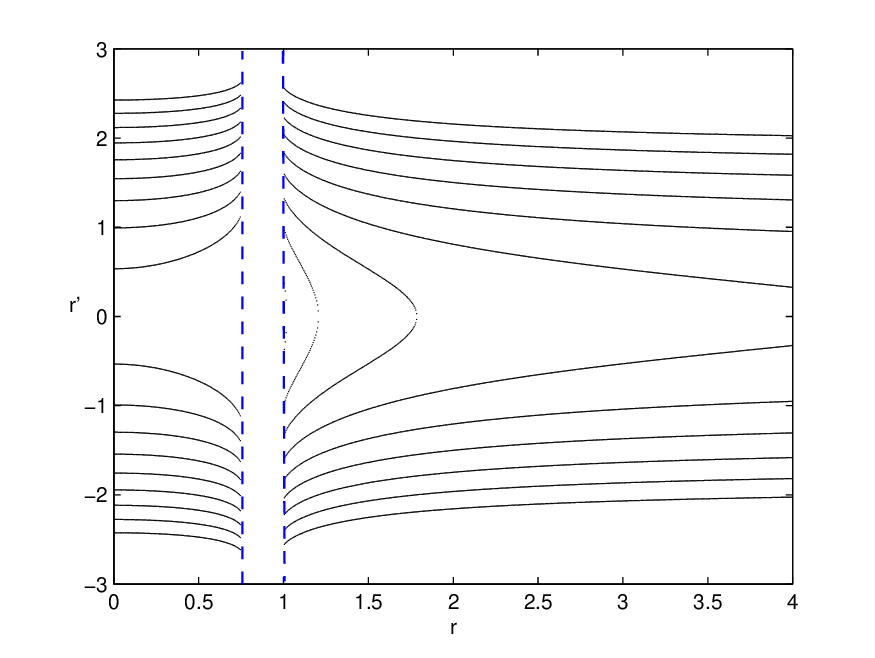}
  \caption{Representation of $W'(r)$ and phase portrait when $\Lambda=0$}\label{fig:derWeff_r_c0}
\end{figure}

 In order to find stationary points we begin with by plotting and analyzing some figures. Fig.~\ref{fig:derWeff_r_c0} (left) depicts the derivative of the effective potential \Eqref{eqMotion} when $\Lambda=0$. We observe that the origin $r=0$ is a critical point, and so, the dynamics are reduced to the linear motion along a diameter of the annulus. There is another zero inside the annular disk ($b<r<a$); although it may appear that orbits crossing the disk have no physical meaning, let us recall that in fact, planetary rings are made of millions of rocky and icy particles, each maintaining their own orbit around the planet, these small orbiting particles can be considered, from a distance, as a continuous solid annular ring. In those cases the critical point would have a physical meaning, however these orbits will end up quickly in a collision orbit.

\begin{figure}[h]
\centering
\includegraphics[width= 2.3in]{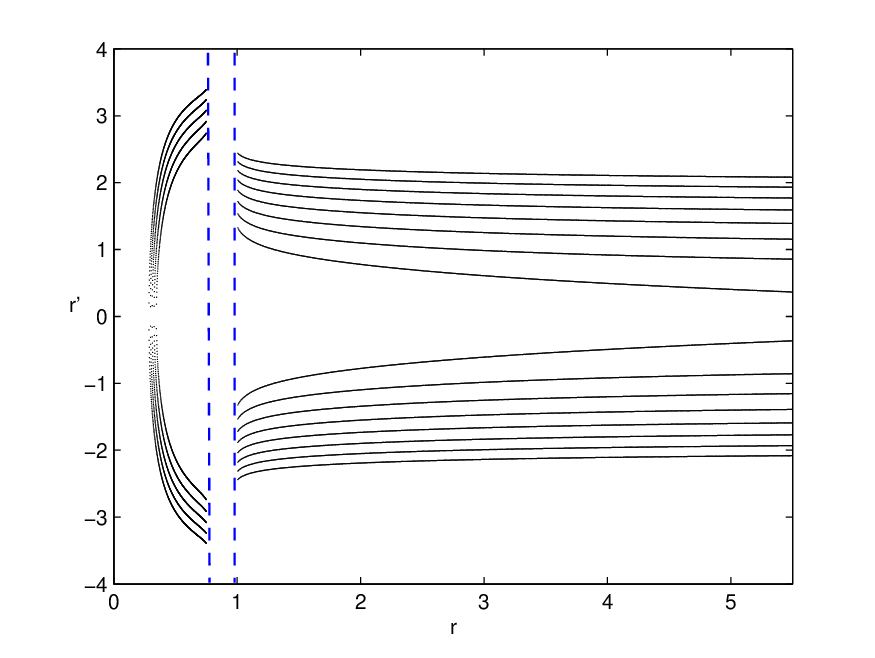}
\includegraphics[width= 2.3in]{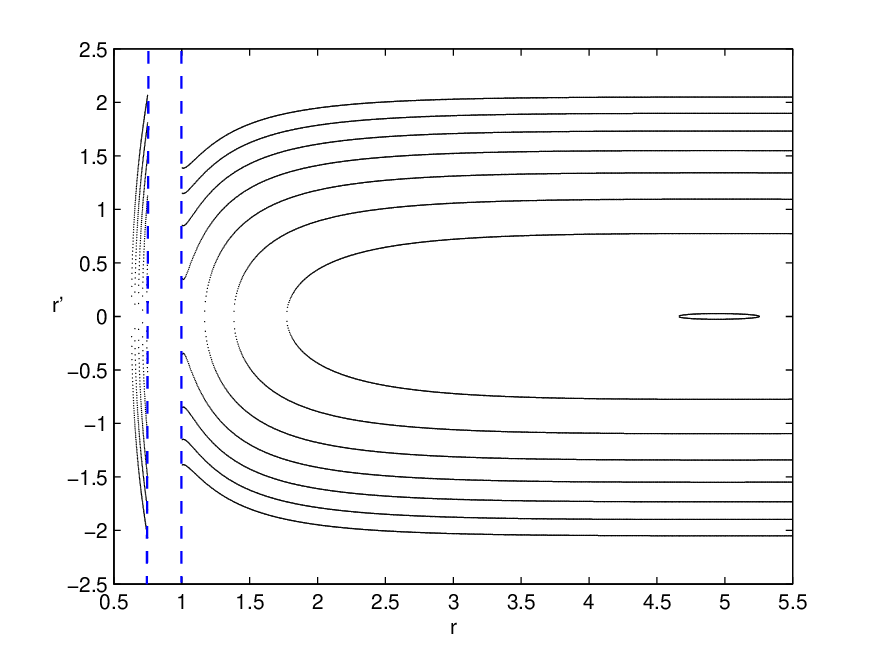}
  \caption{Phase portrait with $\Lambda=1$ and $\Lambda=2.25$ respectively}\label{fig:phase_port_LA}
\end{figure}

This is also observed in the phase portrait $(r,\dot{r})$, Fig.~\ref{fig:derWeff_r_c0} (right); note that inside the annulus we find orbits that originate and end colliding with the annulus, while outside we also find scape orbits. Nevertheless, the phase portrait when $\Lambda\neq 0$ seems to indicate different dynamics depending on the value of the angular momentum, see Fig.~\ref{fig:phase_port_LA}.

\begin{figure}[h]
\centering
\includegraphics[width= 3in]{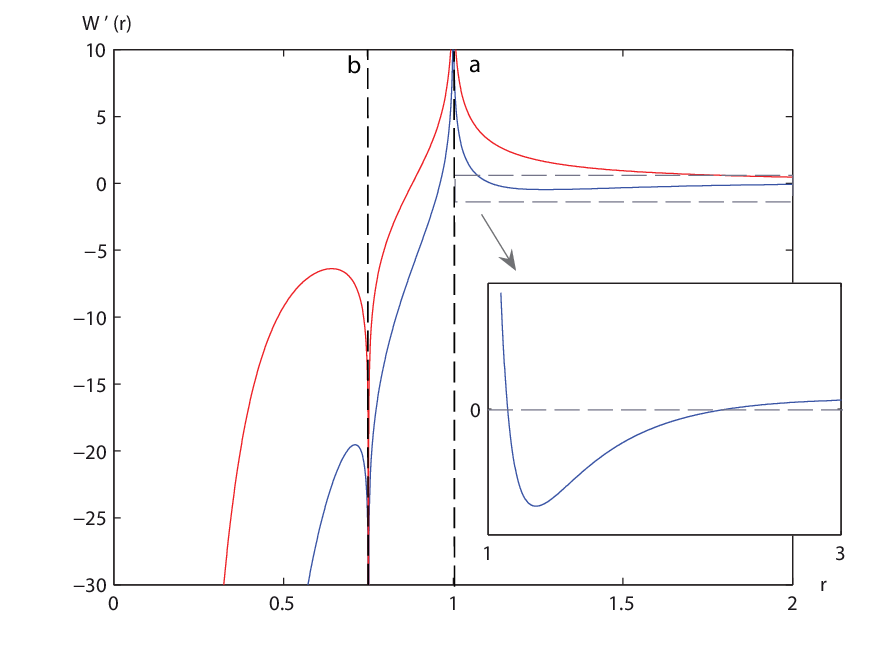}%cor_derWeff_r_mns
  \caption{Critical points of the effective potential for two different values of $\Lambda\neq 0$.}\label{fig:derWeff_r}
\end{figure}

Outside the annulus, and depending on the angular momentum we will have either similar behavior than in the case $\Lambda=0$, or we find a critical point for increasing values of $r$ and another one quite close to the annulus. If we plot the derivative of the effective potential for different values of the angular momentum, (Fig.~\ref{fig:derWeff_r}), we can observe that there are no critical points inside the annulus for any value of  $\Lambda$ and a particle placed there will tend to collide. Outside the annulus we see again that depending on the value we get either none or two zeros.

Therefore, we have numerically found that there exists a bifurcation value $\tilde{\Lambda}$ so that for smaller values $\Lambda \leq \tilde{\Lambda}$ there are no critical points, while for $\Lambda>\tilde{\Lambda}$ the effective potential presents a local maximum and a minimum corresponding to circular orbits of the problem.

We  summarize here briefly the main steps of the analytical proof of this fact. We perform an analysis of the sign and monotony properties of the effective potential in the inner and outer regions of the annulus. Instead of using the expression given in Eq.~\Eqref{eq:cilEqMot_r_z0} we will make use of an integral expression of the potential $U(r)$ when the point $P$ where we compute the potential is on the $xy$-plane
\[
U(r) = -\frac{2\mu}{\pi(a^2-b^2)}\int_{b}^{a}\int_{0}^{\pi}\frac{\rho d\vartheta d\rho} {\sqrt{r^2+\rho^2-2\rho r\cos2\vartheta}}.
\]
Here $r=\sqrt{x^2+y^2}$ is the distance from $P$ to the origin of coordinates, the arbitrary point over the annulus $Q$ is at distance $\rho$ from the origin, and $\vartheta$ is the angle $\widehat{QOP}$. Now changing the variable of integration $\theta=\vartheta+\pi/2$ and after some trigonometric transformations, the potential can be written as
\[
U(r) = -\frac{4\mu}{\pi(a^2-b^2)}\int_{b}^{a}\frac{\rho}{r+\rho} K(k)d\rho,\]
where $k^2= 4\rho r/(r+\rho)^2$. Hence
\[
U'(r) = -\frac{4\mu}{\pi(a^2-b^2)}\int_{b}^{a}\left(\frac{-\rho}{(r+\rho)^2} K(k) + \frac{\rho}{(r+\rho)}\frac{4\rho(\rho-r)}{(r+\rho)^3}\frac{dK(k)}{dk}\right)d\rho.
\]
To determine the sign of the function $U'(r)$ it is enough to study sign of the integrand function $f$, if positive, it means that its primitive is an increasing function in the interval $b<a$ so the integral will be positive; and respectively for the negative sign. Thus, taking into account the properties of the Elliptic Integrals
\[
K(t)>0,\quad \frac{dK(t)}{dt} = \frac{E(t)-(1-t)K(t)}{2t(1-t)}>0,\quad\forall 0<t<1
\]
we can prove that the derivative of the potential function, $U'(r)$, is a negative function in the inner region of the annulus $r<b$, and positive in the outside it $r>a$.

The proof is clear in the case $r>a$, taking into account that $\rho\in[b,a]$ and the properties written above that guarantee $U'(r)>0$. To prove the case $r<b$ we need to replace the derivative of the Elliptic Integral of first kind, and so the derivative of the potential yields
\[
U'(r) = -\frac{4\mu}{\pi(a^2-b^2)}\int_{b}^{a}\left( \frac{-\rho}{2r(r^2-\rho^2)}\big[(r+\rho)E(k)+(r-\rho)K(k)\big]\right)d\rho.
\]
Finally note that the expression $(r+\rho)E(k)+(r-\rho)K(k)$ is the module of the normalized potential Eq.~\Eqref{eq:potefe} at points on the Equatorial plane, and so, it can not be negative. Thus the integrand is a positive function and we can conclude that inside the annulus $U'(r)<0$.
\par

This result allows us to state straightforwardly that there are no critical points of the effective potential inside the annulus, which was already observed in Figs.~\ref{fig:derWeff_r_c0} and \ref{fig:derWeff_r}. As it was stated above, outside the annulus we can find critical points depending on the value of the angular momentum; we have computed analytically a minimum value of $\Lambda$ that guarantees the existence of zeros, let us see how to proceed.

First, we prove the following properties:
\begin{itemize}
  \item The effective potential is an increasing function when $r$ is close to the outer radius $a$.
  \item $W(r)$ is  negative and increasing function when $r$ tends to infinity.
\end{itemize}

The proof in the first case is straightforward and follows from the derivative of the potential Eq.~\Eqref{eq:cilEqMot_r_z0}, calculating the limit when $r$ goes to $a$. For the second property just note that outside the annulus, according to the previous results about the potential, the following inequalities are satisfied
\[
W'(r)>-\frac{\Lambda^2}{r^3},\quad W(r)<\frac{\Lambda^2}{2r^2},
\]
and so $W(r) < 0$ and $W'(r) > 0$ when $r$ tends to infinity.

 These properties establish that outside the annulus, the effective potential starts increasing for values of $r$ close to $a$, and as the radius $r$ grows the potential goes to zero with negative values. Which means that if we find a value of the angular momentum such that $W(r)>0$ for any $r>a$, we will ensure that the effective potential has at least a local minimum and a local maximum. We can compute a rough bound value for the angular momentum satisfying this, computing the limit of $W(r)$ when $r$ tends to $a$
\[W(r)=\frac{\Lambda^2}{r^3}+\frac{2\mu}{\pi(a^2-b^2)}\Big((b+r)E\left(k_b\right)+(b-r)K\left(k_b\right)
-(a+r)E\left(k_a\right)-(a-r)K\left(k_a\right)\Big).
\]
Taking into account that $\displaystyle \lim_{x\to 1}E(x)=1$ and $\displaystyle \lim_{x\to 1}(1-x)K(x)=0,$ %Eq.~\Eqref{eq:g} is a positive function, and
it follows that
\[
 \lim_{r \to a^+}W(r)=\frac{\Lambda^2}{a^3}+\frac{2\mu}{\pi(a^2-b^2)}\left(-2a+(a+b)E\left(\frac{4ab}{(a+b)^2}\right)+(b-a)K\left(\frac{4ab}{(a+b)^2}\right)\right).
\]
Hence,
\[
\lim_{r \to a}W(r)>0\quad \Leftrightarrow \quad\frac{\Lambda^2}{a^3}-\frac{4\mu a}{\pi(a^2-b^2)}>0\quad \Leftrightarrow \quad \Lambda^2>\frac{8\mu a^3}{\pi(a^2-b^2)}.
\]

For greater values of the angular momentum, we can ensure that we have found a positive value of the effective potential, and so a local minimum and maximum of it. We can summarize all the results obtained in the case of movement reduce to the $xy$-plane in the following proposition.\\

\begin{proposition} The dynamics of a particle under the attraction of a massive annulus and confined on the Equatorial plane verifies that
\begin{itemize}
  \item[a)]  The origin is only one critical point of the effective potential $W(r)$ when the angular momentum $\Lambda=0$.
  \item[b)]  There are no circular solutions in the inner region of the annular disk.
  \item[c)]  For values of the angular momentum $\Lambda^2 > \tilde{\Lambda}^2=8\mu a^3/\pi(a^2-b^2)$ the effective potential has at least two critical points in the outer region of the annulus, corresponding to one stable and one unstable circular orbits.
  \item[d)]  There is a stable critical point inside the annulus $b<r<a$ for any value of the angular momentum.
\end{itemize}
\end{proposition}

We have also observed that as the angular momentum increases, one of the critical points tends to the annulus while the other goes to infinite.

\section{Composition of annulus}
In this section we analyze the dynamics created by two concentric annuli and how the second annulus modifies the results found in the previous section.

We consider two concentric annuli, that is, a partition of the annulus in two smaller ones with a gap between them. This approach provides a better approximation to real planetary rings, and so the study of the dynamical structure around this model is considered of special interest.
The potential function is obtained by adding to the potential created by an annulus of radius $0<b<a,$ given by Eq.~\Eqref{pot_ann}, a second annulus of equal mass and new radius $d<c$.

We have performed an analogous procedure to the one above described to determine the critical points of the effective potential for this model. Since the potential function is an additive function, the qualitative behavior in the inner part or in the outer part of both annuli will be the same as described in the previous section; thus, we focus our interest in the the gap between both annuli.

The partition of the original annulus in two smaller ones implies the existence of a new equilibrium point located inside the gap. The proof is straightforward by computing the limits of the effective potential when the particle approaches the boundary of the annular rings, and taking into account that the derivative of the potential of a single annulus Eq.~\Eqref{eq:cilEqMot_r_z0} goes to infinity when $r$ goes to the outer radius of the annulus, and goes to $-\infty$ when $r$ approaches the inner radius. Thus,  if $W(r)$ denotes the effective potential of the composition of annulus, there results that
\begin{equation}
\begin{array}{l}
\displaystyle\lim_{r\to b^-}W'(r) = -\Lambda^2/b^3 + k/b^2 -\infty = -\infty,\\[2ex]
\displaystyle\lim_{r\to c^+}W'(r) = -\Lambda^2/c^3 + k/c^2 + \infty = \infty.
\end{array}
\end{equation}

\begin{figure}[h]
\includegraphics[width= 10cm]{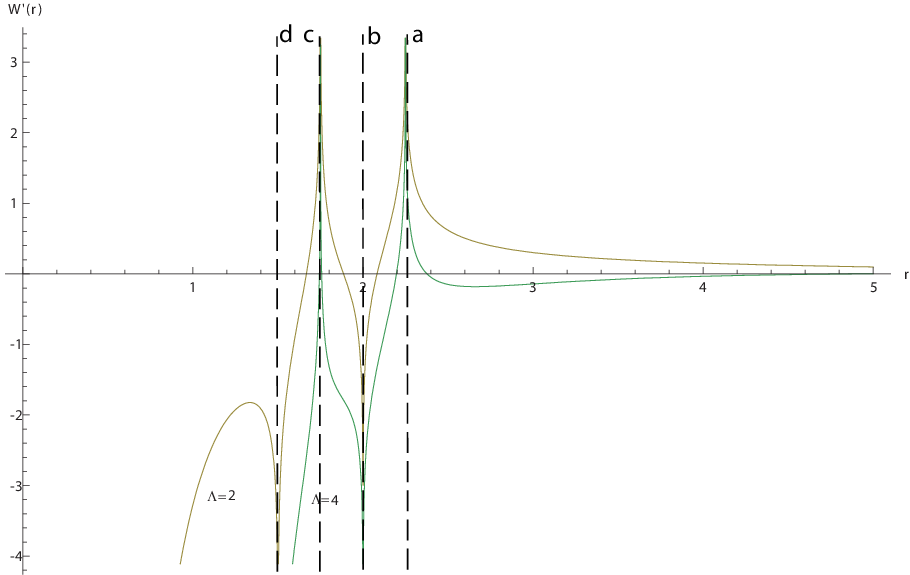}
  \caption{Critical points of the effective potential of two concentric annulus for different values of $\Lambda\neq0$. The yellow line is for $\Lambda = 2$, and the green one is for $\Lambda = 4$.}
  \label{fig:derWeff_r_2}
\end{figure}

The effective potential is a continuous function in the interval between the annular disks (c,b), hence we can conclude that there is at least a critical point (an odd number of equilibria as a matter of fact) in the region confined between both annulus. Besides, by plotting the derivative effective potential (Fig. \ref{fig:derWeff_r_2}) we discover that there is only one zero indeed. This equilibrium point correspond to one circular unstable orbit. The stability has been numerically determined by means of the monodromy matrix, since there is no feasible to do it analytically.

Since this critical point is spectrally unstable, joined to the fact of the existence of a stable equilibrium inside each annulus, it may explain how the dynamics of particles inside and around  these annular systems is organized.

 Figure \ref{fig:derWeff_r_2} represents the derivative of the effective potential for different values of the angular momentum.

\section*{Conclusion}
We have set up the problem of the motion of a massless particle around a massive annulus by finding a closed expression of the potential function valid for all points in space where it is defined. This closed expression, despite its complex treatment due to the Elliptic Integrals involved, avoid approximations of the potential function through asymptotic expansions.

We have also analyzed the dynamics of an infinitesimal particle under the attraction of a planar annulus and a composition of annuli performing a systematic search of one of the most relevant solutions: equilibrium positions. We have proved and given conditions for the existence of equilibrium points for different values of the angular momentum on the Equatorial plane and the polar axis. We have also observed different regions of planar stable motion around the ring, and their possible implication in the formation of dynamical patterns.

\section*{Acknowledgments}
We appreciate comments and suggestions from Dr. Kondratyev who improved the clarity and
structure of the paper, as well as for providing us references of his previous works on the topic. This paper has been supported by the Spanish Ministry of Science and Innovation, Project \# AYA2008-05572.

%%%%%%%%%%%%%%%%%%%%%%%%%%%%%%%%%%%%%%%%%%%%
%%%%%%%%%%%%%%%%%%%%%%%%%%%%%%%%%%%%%%%%%%%%

\end{document}